\documentclass[11pt]{article}
\usepackage{amsmath,amsfonts,amssymb, amsthm}
\usepackage[english,french]{babel}

\newtheorem{theorem}{Theorem}
\newtheorem*{maintheorem}{Main Theorem}
\newtheorem*{theorem*}{Theorem}
\newtheorem{proposition}{Proposition}
\newtheorem*{proposition*}{Proposition}
\newtheorem{lemma}{Lemma}
\newtheorem{remark}{Remark}
\newtheorem{corollary}{Corollary}

\newcommand{\R}{\mathbb{R}}

\newcommand{\C}{\mathbb{C}}
\newcommand{\Q}{\mathbb{Q}}
\newcommand{\D}{\mathbb{D}}
\newcommand{\Fc}{\mathcal{F}}

\newcommand{\eps}{\varepsilon}
\newcommand{\La}{\mathcal{L}}

\newcommand{\PP}{\mathbb{P}}

\title{A Characterization of Harmonic Measures on Laminations by Hyperbolic Riemann Surfaces}
\author{Yuri Bakhtin\thanks{Georgia Institute of Technology and Fields Institute} 
\and Matilde Mart\'\i nez\thanks{Fields Institute and CIMAT}}

\begin{document}

\selectlanguage{english}

\maketitle

\begin{abstract} $\La$ denotes a (compact, non-singular) lamination by hyperbolic Riemann surfaces. 
We prove that a probability measure
on $\La$ is harmonic if and only if it is the projection of a measure
on the unit tangent bundle $T^1\La$ of $\La$ which is invariant under both
the geodesic and the horocycle flows.
\end{abstract}
\begin{otherlanguage}{french}
\begin{abstract}
$\La$ denote une lamination (compacte, non singuli\`ere) par surfaces de Riemann hyperboliques. On montre qu' une mesure sur $\La$ est harmonique si et seulement si elle est la projection d'une mesure sur le fibr\'e tangent unitaire $T^1\La$ qui est invariante sous les flots g\'eodesique et horocyclique.
\end{abstract}
\end{otherlanguage}


\section*{Introduction}

This note studies measures associated to compact nonsingular laminations by hyperbolic Riemann surfaces. 
It builds on a previous article by the second author, entitled~\emph{Measures on Hyperbolic Surface Laminations} 
(see~\cite{Martinez}), and improves on the main results found therein. 

In this work,~$\La$ denotes a lamination (or foliated space), which is compact and whose leaves 
are hyperbolic Riemann surfaces. Each leaf has its Poincar\'e metric -- the metric of constant curvature $-1$
compatible with the conformal structure. We consider two different kinds of measures associated to~$\La$: 
measures invariant under the heat diffusion along the leaves of $\La$ which are called \emph{harmonic} 
and measures on the unit tangent bundle~$T^1\La$ of~$\La$ that are invariant under the laminated geodesic 
and horocycle flows. Our main result reads as follows:

\begin{maintheorem} Let $\mu$ be a probability measure on $\La$. Then $\mu$ is harmonic if and only if there is a measure  $\nu$ on $T^1\La$ which is invariant under both the geodesic and stable horocycle flow and that projects onto $\mu$ (under the canonical projection $T^1\La\rightarrow\La$). Furthermore, such a $\nu$ is unique.
\end{maintheorem}

This had been conjectured by Christian Bonatti.

Of course, in the statement of this theorem we can change the \emph{stable} horocycle flow for the \emph{unstable} horocycle flow, but asking that the measure $\mu$ be invariant under the three flows we are considering might be too much: In fact, there is a canonical bijection between measures invariant under the three flows and the simplest kind of harmonic measures; namely, those coming from measures on transversals which are invariant under the holonomy pseudogroup. This is easy to prove and can be found in \cite{Martinez}.

Our Main Theorem represents an improvement on the main results found in \cite{Martinez} (labeled Theorem 2.6 and Theorem 2.8), which state the following:

\begin{theorem*}
\begin{enumerate}
\item  Any harmonic probability measure on $\La$ is the projection of a measure invariant under the stable horocycle flow on $T^1\La$.

\item  Any probability measure on $T^1\La$ invariant under both the geodesic and the stable horocycle flow projects onto a harmonic measure on $\La$.  
\end{enumerate}
\end{theorem*}

Remark that the first statement is a partial converse of the second one; what we do in this note is to prove the converse of 2. The techniques we use are completely different from those used in \cite{Martinez}.

In Section 1 we review and slightly modify the construction of a measure on $T^1\La$ invariant 
under the horocycle flow that projects onto a given harmonic measure on~$\La$. 
We wish to prove that the measure thus constructed is also invariant under the geodesic flow.  
This can be reduced to a probabilistic statement on the regularity
of the distribution of the radial component of the Brownian motion on the hyperbolic plane.

The key regularity statement is provided by Theorem~2 which is proved in Section~2. The main idea
of the proof is to show that for large times the radial component of the hyperbolic planar Brownian
motion behaves as one-dimensional Brownian motion with constant positive drift. The closeness estimates
are based on Girsanov's theorem.

The uniqueness of the measure $\nu$ is proved in section 3. We consider a flow box in which we disintegrate a harmonic measure $\mu$ using Rokhlin's theorem, and prove that there is at most one way to lift the harmonic measure in this flow box to a measure invariant under the affine group. Only local considerations are involved.

In section 4 we state some simple applications of our Main Theorem. We state previous results regarding the uniqueness of harmonic measures which can now be translated to unique ergodicity for the action of the affine group in some examples. 

Acknowledgment. The authors gratefully acknowledge the hospitality of the Fields Institute in Toronto where this work has been done.
The second author is grateful for partial support from NSERC, via the research grant of M.~Lyubich.
We would also like to thank M.~Lyubich personally for his input and comments. We are deeply grateful to both referees for their corrections, comments, and a very careful reading of the manuscript.

\section{Harmonic measures as projections of measures invariant under the affine group.}

\subsection{Hyperbolic surface laminations}

$\La $ is a \emph{lamination} if it is a separable,
locally compact
metrizable space that has an open covering $\{ E_i\}$ and an
atlas $\{(E_i,\varphi_i)\}$ satisfying:

\begin{enumerate}
\item $\varphi_i:E_i\rightarrow D_i\times T_i$ is a homeomorphism, for some open disk $D_i$ in $\R^d$ and topological
space $T_i$, and

\item the coordinate changes $\varphi_j\circ\varphi_i^{-1}$ are  of the form
$(z,t)\mapsto (\zeta (z,t),\tau(t))$, where each $\zeta$ is smooth in the
$z$ variable. 
\end{enumerate}

This last
condition says that the sets of the form $\varphi_i^{-1}(D_i\times \{t\})$, called
\emph{plaques}, glue together to form $d$-dimensional manifolds that we call
\emph{leaves}.

$\La $ is a \emph{Riemann surface lamination}
if the disks $D_i$ are open subsets
of the complex plane and the maps $\zeta$ are holomorphic in the $z$ variable.
We say that $\La$ is a \emph{hyperbolic surface lamination} if its
leaves are hyperbolic Riemann surfaces.

Each leaf on a hyperbolic surface lamination $\La$ has a Poincar\'e metric, which is
the only Riemannian metric of constant curvature $-1$ compatible with the conformal
structure. According to a theorem
due to Candel (see \cite{Candel1}), these metrics on the leaves, as well as all their derivatives,
have
continuous variation in the transverse direction.

In this paper, $\La $ will always denote a compact hyperbolic surface lamination.

\subsection{Harmonic measures}  

Each leaf $L$ of $\La$, being a Riemannian manifold, has a  Laplace-Beltrami operator $\Delta_L$.  If
$f:\La\rightarrow\R$ is a function of class $C^2$ in the
leaf direction and $x\in\La$, we define $ \Delta
f(x)=\Delta_L f|_L(x),$ where $L$ is the leaf passing through $x$
and $f|_L$ is the restriction of $f$ to $L$.  A probability measure
$\mu$ on $\La$ is \emph{harmonic} if
$\Delta \mu=0$; i.e. if $\int \Delta f\, d\mu=0\ \forall f$. 

The operator $\Delta/2$ is the infinitesimal generator of the \emph{laminated heat semigroup}
$\{ D_t\} _{t\geq 0}$.
(This choice of the generator corresponds to the
standard hyperbolic planar Brownian motion with unit diffusion.)
The  heat diffusion along the leaves can be described in terms of the \emph{laminated heat kernel}, $p:\La\times
\La\times (0,+\infty )\rightarrow \R$ given by
\begin{equation*}
p(x,y,t)=
\begin{cases}p_{L}(x,y,t),& \mbox{if}\ x,y\ \mbox{belong to the same leaf}\ L\\ 0,& \mbox{otherwise,}
\end{cases}
\end{equation*}
where $p_L$ is the heat kernel on $L$.
The laminated heat semigroup  can be expressed as
$$D_tf(x)=\int_{L_x} p(x,y,t)f(y)\, dy,$$
where $L_x$ is the leaf through $x$.

Lucy Garnett proved an Ergodic Theorem for harmonic measures, which can be found in \cite{Garnett} or \cite{Candel2}.

A harmonic measure
$m$ is \emph{ergodic} if $\La$ can not be partitioned into two
measurable leaf-saturated subsets having positive $m$ measure.

The ergodic harmonic measures can be briefly described as follows: we take a
point $x$ in $\La$ and consider the Dirac delta at $x$, which we call $\delta _x$.
We can diffuse this measure, obtaining for each positive time $t$ a probability
measure $D_t\delta_x$, whose integral on any continuous function $f$ on $\La$ is
$D_tf(x)$. For almost all $x$ according to any harmonic measure, the sequence
of Krylov--Bogolyubov means
$$\frac{1}{n}\sum_{t=0}^{n-1}D_t\delta_x$$
has a limit which we call $\tilde{\delta_x}.$ Ergodic measures are 
of the form $\tilde{\delta_x}$.

\subsection {The $PSL(2,\R)$-action on the unit tangent bundle}

If $\La$ is a hyperbolic surface lamination, we call $T^1\La$ the lamination whose
three dimensional leaves are the unit tangent bundles of the leaves of $\La$ and
that has ``the same" charts as $\La$.

The \emph{laminated geodesic flow} is the flow $g_t$ that,
restricted to the unit tangent bundle of a leaf
$L$ of $\La$, coincides with the geodesic flow in $T^1L$. The \emph{laminated
stable and unstable horocycle flows} $h^+$  and $h^-$ are the
flows that, when restricted to the unit tangent bundle of a leaf $L$, coincide with
the stable and unstable horocycle flows on $T^1L$, respectively.
For the definition and basic properties of the geodesic, stable horocycle
and unstable horocycle flows on hyperbolic surfaces, see \cite{Manning}.
All these flows are
continuous on $T^1\La$, as a consequence of Candel's theorem (see \cite{Candel1}).

In the lamination $T^1\La $ there is a
right $PSL(2,\R)$-action
whose orbits are the leaves. As in the case of surfaces, the geodesic and the
horocycle flows correspond to the action on $T^1\La $ of the one-parameter subgroups $D_t$, $H^+$ and $H^-$, respectively, where 

$$D=\left\{\left(\begin{array}{cc}e^{\frac{t}{2}}&0\\ 0& e^{-\frac{t}{2}}
\end{array}\right)\right\},\ \ 
H^+=\left\{\left(\begin{array}{cc} 1&t\\ 0&1
\end{array}\right)\right\},\ \  
\hbox{ and }\ 
H^-=\left\{\left(\begin{array}{cc} 1&0\\ t&1
\end{array}\right)\right\}.$$

Consequently, the joint action of $D$ and $H^+$ corresponds to the action of the affine group

$$B=\left\{\left(\begin{array}{cc} a&b\\ 0& a^{-1}
\end{array}\right);\ a,b\in\R, a> 0\right\}.$$

\subsection{A measure invariant under the action of the affine group}

Let $\pi:T^1\La\rightarrow\La$ be the canonical projection, and consider a harmonic probability measure $\mu$ on 
$\La$. In \cite{Martinez}, there is a construction that produces a measure $\nu$ on $T^1\La$ which is invariant under the 
stable (or the unstable) horocycle flow and such that $\pi_*\nu=\mu$.
It will be more convenient for us to consider the measure $\nu$ which is invariant 
under the \emph{unstable} horocycle
flow and which projects onto $\mu$.

\begin{theorem} The measure $\nu$ is invariant under the geodesic flow $g$. \end{theorem}

This theorem says that $\nu$ is, in fact, invariant under both the geodesic and horocycle flows. 
Therefore, any
harmonic measure is the projection of a measure invariant under both flows. This
statement and its converse which is quoted in the Introduction constitute our Main Theorem.

{\sc Proof of Theorem 1.}
As in \cite{Martinez}, two simplifying assumptions shall be made, that imply no loss of generality:

\begin{enumerate}
\item That $\mu$ is ergodic; that is, there is an $x\in\La$ such that $\mu=\tilde\delta_x$; and

\item that the point $x$ belongs to a leaf which is simply connected. (If it does not, we consider the universal cover of the leaf through $x$, as explained in \cite[pg.857]{Martinez}.)
\end{enumerate}

The construction of $\nu$ goes as follows:

For any natural number $n\geq 1$, let $\delta^{(n)}_x$ be the Krylov--Bogolyubov sum $\frac{1}{n}\sum _{t=0}^{n-1} D_t\delta_x$. Namely, $\delta^{(n)}_x$ is the probability 
measure such that, for every continuous function $f$ in $\La$,
$$\int f\, d\delta^{(n)}_x =\frac{1}{n}\sum_{t=0}^{n-1} D_tf(x)=\frac{1}{n}\sum_{t=0}^{n-1}\int _\La p(x,y,t)f(y)dy.$$
With this notation, $\mu=\lim_n \delta^{(n)}_x$.

Let $L_x$ be the leaf of the lamination $\La$ passing through $x$. It is a hyperbolic plane. Let 
$R:L_x\backslash \{x\}\rightarrow T^1\La$ be the unit radial vector field pointing outwards; i.e. 
$R(y)=(\gamma(s),\dot\gamma(s))$ if $\gamma$ is the geodesic of unit speed such that $\gamma(0)=x$ and 
$\gamma(s)=y$. (Here $s$ is, of course, the distance from $x$ to $y$ measured on $L_x$.)

Define
$\mu_n=R_*\delta^{(n)}_x$. This gives a sequence of probability measures on the compact space $T^1\La$. Let 
$\nu$ be any limit point of the sequence $(\mu_n)$ in the sense of the weak-* topology. (It is not difficult to 
see that in fact $\nu$ is the limit of the $\mu_n$, but this is unessential for our argument.) Then, as was proved in 
\cite{Martinez}, $\nu$ is invariant under the unstable horocycle flow. And $\pi_*\nu=\mu$ since $\pi_*$ is a continuous map 
from the space of finite measures on $T^1\La$ to that of finite measures on $\La$.

To finish the proof of Theorem 1 we shall show that $(g_s)_*\nu=\nu$ for all $s\in [0,1]$. It is clearly enough to prove the following:

(*) For every continuous real-valued function $f$ on $T^1\La$ and every $s\in [0,1]$,
$$\lim_{t\rightarrow +\infty}\left|\int_{T^1\La}f\, d(R_*D_t\delta_x)-\int_{T^1\La}f\circ g_s\, d(R_*D_t\delta_x)\right|=0.$$

Let $S_r=S_r(x)=\{ y\in L_x :\, d(x,y)=r\}$,
where $d$ is the hyperbolic distance measured on the leaf $L_x$. Its normalized Lebesgue measure can be pushed forward 
by $R$ to get a measure, that we call $\lambda_r$, supported on the curve $R(S_r)\subset T^1\La$.
We can write the integral of $f$ with respect to $R_*D_t\delta_x$ as
$$\int_{\R^+} p_t(r)\left(\int f\, d\lambda_r \right) \, dr,$$
where 
$p_t(r)=p(x,y,t)\times (\hbox{length}(S_r))$ for any point $y$ in $L_x$ such that $d(x,y)=r$.
If $s>0$, the geodesic flow at time $s$ takes $R(S_r)$ to $R(S_{r+s})$, and $(g_s)_*\lambda_r=\lambda_{r+s}$.
Writing $u(r)=\int f\, d\lambda_r$, the statement labeled (*) follows from:

(**) For every continuous bounded function $u:\R^+\rightarrow\R$ and every $s\in [0,1]$,
$$\lim_{t\rightarrow\infty}\left|\int_{\R^+} p_t(r)u(r)\,dr - \int _{\R^+} p_t(r+s)u(r)\, dr\right|=0.$$

Claim (**) is implied by

\begin{theorem} For any $s\in [0,1]$,
\begin{equation}
\label{eq:closeness_of_densities}
\lim_{t\rightarrow\infty}\int_{\R^+}|p_t(r)-p_t(r+s)|dr=0.
\end{equation} 
\end{theorem}

Section 2 is devoted to the proof of this theorem. $\square$

\begin{remark} Theorem 1 says that any harmonic measure on $\La$ is the 
projection of a measure invariant under the action of the lower triangular group in $T^1\La$. Of course
we could do the same thing for the upper triangular group, considering the unit inward radial vector field
on $L_x\backslash \{x\}$, and then taking $s$ in $[-1,0]$.\end{remark}


\section{The regularity of the heat kernel} 
Though an explicit expression for $p_t(r)$ is well-known, see e.g.~\cite[p.246]{Chavel},
it is not easy to use that formula directly. We use a stochastic calculus approach instead.
This method can be generalized to higher  dimensions.

We start with the following rotationally invariant Riemannian metric on the plane $\R^2$
in polar coordinates:
\begin{equation*}
ds^2=dr^2+\sinh^2(r)d\theta^2
\end{equation*}
which turns it into a hyperbolic plane (see \cite[Chapter X]{Chavel}).

In our study of the Brownian motion on the hyperbolic plane  we shall refer to
\cite{Karatzas-Shreve} and \cite{RevuzYor} for basic facts on Brownian motion and stochastic
calculus.

Let $B_t$ be the standard Brownian motion on this hyperbolic plane (a stochastic process with
generator given by $\Delta/2$ where $\Delta$ is the Laplace--Beltrami operator associated with $ds^2$).
Its radial component $X_t=d(0,B_t)$ satisfies the following stochastic
differential equation:
\begin{equation}
\label{eq:main-SDE}
dX_t=dW_t+\frac{1}{2}\coth (X_t)dt
\end{equation}
where $W_t$ is a standard Wiener process. This can be easily derived using the explicit
expression for the Laplace--Beltrami operator in polar coordinates $(r,\theta)$, see e.g.~\cite{March}.
Equation \eqref{eq:main-SDE} should be understood in the integral sense:
\begin{equation}
\label{eq:main-SDE-integral}
X_t-X_0=W_t+\frac{1}{2}\int_0^t\coth(X_s)ds.
\end{equation}
This equation defines a family of transition probability kernels:
\begin{equation}
P_{t}(x,dy)=\PP\{X_t\in dy|\ X_0=x\},\quad t>0, x\in\R_+.
\end{equation}
If $x\in\R_+$, then $P_{t}(x,\cdot)=1$ is concentrated on $\R_+$
due to a singularity of the drift term at $0$. (This singularity point
is of ``entrance and non-exit'' type according to the classification in
~\cite{Ito-McKean}).

These kernels are, in fact, absolutely continuous with respect to 
the Lebesgue measure $dy$ on $\R_+$, and the associated density $p_t(x,y),x,y\in\R_+$ satisfies
the following forward Kolmogorov (or Fokker--Planck) equation (see~\cite[equations~(1.6) and (1.8) on p.282]{Karatzas-Shreve}):
\begin{equation}
\frac{\partial p_t(x,y)}{\partial t}=\frac12 \frac{\partial^2 p_t(x,y)}{\partial y^2}-
\frac{\partial}{\partial y}\left(\frac{1}{2}\coth(y) p_t(x,y)\right), \quad x,y\in\R_+.
\end{equation}

For $y\in\R_+$ we shall denote $p_t(y)=p_t(0,y)$ which is consistent with the definition
of $p_t(\cdot)$ given in Section 1.4, since the heat kernel $p(0,\cdot,t)$ is rotationally invariant and 
coincides with the
transition density of the hyperbolic Brownian motion.

\begin{remark} In probabilistic terms, equation~\eqref{eq:closeness_of_densities} from the statement
of Theorem~2 is equivalent to
\begin{equation*}
\lim_{t\to\infty}\|P_t(\cdot)-P_t(\cdot+s)\|_{TV}=0,
\end{equation*}
where $P_t(\cdot)=P_t(0,\cdot)$, $P_t(\cdot+s)$ means $P_t(\cdot)$ translated by $s$,
and $\|\mu_1-\mu_2\|_{TV}$ denotes the total variation distance between measures $\mu_1$
and $\mu_2$, see \cite[Section 5.3]{Thorisson}.
\end{remark}

{\sc Proof of Theorem 2.}
Our strategy will be to compare  $p_t(y,z)$ to the transition density of
the Brownian motion with constant drift $1/2$. This will be possible due to
Girsanov's theorem and closeness of the drift term $\coth(x)/2$ in equation~\eqref{eq:main-SDE}
to $1/2$ for large values of $x$.

We shall need the following Kolmogorov--Chapman equation (see~\cite[Chapter III]{RevuzYor}):
\begin{equation}
\label{eq:Kolmogorov-Chapman}
P_t(x,A)=\int_{\R_+} P_{t'}(x,dy)P_{t-t'}(y,A),\quad t'\in(0,t),
\end{equation}
or, equivalently,
\begin{equation}
\label{eq:Kolmogorov-Chapman-density}
p_t(x,z)=\int_{\R_+} P_{t'}(x,dy)p_{t-t'}(y,z),\quad t'\in(0,t).
\end{equation}

Let us fix $s\in[0,1]$ and denote the integral under limit on 
the l.h.s. of~\eqref{eq:closeness_of_densities} by $I_t$.
Equation~\eqref{eq:Kolmogorov-Chapman-density} implies
\begin{equation}
\label{eq:estimate1}
I_t\le \int_{\R_+}P_{t'}(dy)\int_{\R_+}|p_{t-t'}(y,z)-p_{t-t'}(y,z+s)| dz,\quad t'\in(0,t).
\end{equation}
 
Next,
\begin{align}
\label{eq:estimate2}
\int_{\R_+}|p_{t-t'}(y,z)-p_{t-t'}(y,z+s)|dz &\le 
\int_{\R_+}|p_{t-t'}(y,z)-p^*_{t-t'}(y,z)|dz
\\&+\int_{\R_+}|p_{t-t'}(y,z+s)-p^*_{t-t'}(y,z+s)|dz \notag
\\&+\int_{\R_+}|p^*_{t-t'}(y,z)-p^*_{t-t'}(y,z+s)|dz \notag
\\&=I_1(y,t-t')+I_2(y,t-t')+I_3(y,t-t'),\notag
\end{align}
where
\begin{equation}
\label{eq:p^*}
p^*_{t}(y,z)=\frac{1}{\sqrt{2\pi t}}e^{-\frac{(z-y-t/2)^2}{2t}}
\end{equation}
is the density of $Y_t=y+W_t+t/2$ which is a Wiener process with constant drift $1/2$ started at $y$.
In fact, $Y_t$ is a solution of the following SDE:
\begin{equation}
\label{eq:constant_drift}
dY_t=dW_t+\frac{1}{2}dt.
\end{equation}
Since $s\in[0,1]$, a straightforward estimate based on~\eqref{eq:p^*} implies that for some constant $K$ and
all $y$,
\begin{equation}
\label{eq:I_3}
I_3(y,t-t')\le\frac{K}{\sqrt{t-t'}}.
\end{equation}

Let us estimate $I_1(y,t)$ for large values of $y$. Namely, let us fix a number $R>0$ to be specified
later and assume that $y>R$.

Girsanov's theorem (see~\cite[Chapter VIII]{RevuzYor}) implies that the measures on paths generated by solutions of~\eqref{eq:main-SDE} 
and~\eqref{eq:constant_drift} emitted from the same initial point $y$ 
are absolutely continuous with respect to the Wiener measure on paths $W$ with
\begin{equation*}
\int_0^t\coth^2(y+W_r)dr<\infty
\end{equation*}
(i.e. on paths not crossing $0$),
and the densities are given by
\begin{equation*}
Z_{t}(y+W)=\exp\left\{\int_0^t\frac{\coth(y+W_r)}{2}dW_r-\frac{1}{2}\int_0^t\left(\frac{\coth(y+W_r)}{2}\right)^2dr\right\}
\end{equation*}
and
\begin{equation*}
Z^*_{t}(y+W)=\exp\left\{\int_0^t\frac{1}{2}dW_r-\frac{1}{2}\int_0^t\left(\frac{1}{2}\right)^2dr\right\}=\exp\left\{\frac{1}{2}W_t-\frac{t}{8}\right\},
\end{equation*}
respectively.

It\^o's formula (see~\cite[p.149]{Karatzas-Shreve}) gives
\begin{equation*}
\ln\sinh(y+W_{t})-\ln\sinh(y)=\int_0^t\coth(y+W_r)dW_r+\frac{1}{2}\int_0^t(1-\coth^2(y+W_r))dr
\end{equation*}
which implies
\begin{equation*}
Z_{t}(y+W)=\left(\frac{\sinh(y+W_{t})}{\sinh(y)}\right)^{1/2}\exp\left\{-\frac{t}{4}+\frac{1}{8}\int_0^t\coth^2(y+W_r)dr\right\}.
\end{equation*}
Therefore
\begin{multline*}
\frac{Z_{t}(y+W)}{Z^*_{t}(y+W)}=\left(\frac{\sinh(y+W_{t})}{\sinh(y)e^{W_t}}\right)^{1/2}
\exp\left\{\frac{1}{8}\int_0^t\coth^2(y+W_r)dr-\frac{t}{8}\right\}\\
=\left(\frac{1-e^{-2(y+W_t)}}{1-e^{-2y}}\right)^{1/2}\exp\left\{\frac{1}{8}\int_0^t(\coth^2(y+W_r)-1)dr\right\}
\end{multline*}
for all paths $y+W$ not crossing $0$.

So, if $y\in[R,\infty)$, and the entire path $y+W$ lies in $[R/2,\infty)$, then
\begin{equation}
(1-e^{-R})^{1/2}\le\frac{Z_{t}(y+W)}{Z^*_{t}(y+W)}\le\frac{\exp\left\{\frac{t}{8\sinh^2(R/2)}\right\}}{(1-e^{-2R})^{1/2}},
\label{eq:density_ratio}
\end{equation}
Let us now split the densities $p^*_t(y,z)$ and $p_t(y,z)$ as follows:
\begin{align*}
p^*_t(y,z)&=q^*_t(y,z)+\psi^*_t(y,z),\\
p_t(y,z)&=q_t(y,z)+\psi_t(y,z),
\end{align*}
Here $q^*_t(y,z)$ (respectively, $q_t(y,z)$) denotes the contribution to $p^*_t(y,z)$ 
(respectively, $p_t(y,z)$) from paths $y+W$ connecting $y$ and $z$
and staying within $[R/2,\infty)$.

From~\eqref{eq:density_ratio} we obtain
\begin{equation*}
\frac{1}{C(R,t)}\le\frac{q_{t}(y,z)}{q^*_{t}(y,z)}\le C(R,t),
\end{equation*}
where
\begin{equation*}
C=C(R,t)= \frac{\exp\left\{\frac{t}{8\sinh^2(R/2)}\right\}}{(1-e^{-R})^{1/2}}.
\end{equation*}

Then,
\begin{equation}
\label{eq:decomposeI_1}
I_1(y,t)\le \int_{\R_+}|q_{t}(y,z)-q^*_t(y,z)|dz + \int_{\R_+}\psi_t(y,z)dz + \int_{\R_+}\psi^*_t(y,z)dz
\end{equation}
The first term is bounded by
\begin{align}
&\int_{z\ge R,\ q_{t}(y,z)>q^*_t(y,z)} (Cq^*_t(y,z)-q^*_t(y,z))dz \notag
\\+&
\int_{z\ge R,\ q_{t}(y,z)\le q^*(y,z)} (Cq_t(y,z)-q_t(y,z))dz \notag 
\\
\le& (C-1)\int_{z\ge R}q_t(y,z)dz \le C-1. \notag 
\label{eq:I_1-main_contrib}
\end{align}

The last term in~\eqref{eq:decomposeI_1} is bounded by
\begin{equation*}
\PP\left\{\min_r (y+W_r+r/2)< R/2\right\}
\le \PP\left\{\min_r (W_r+r/2)<-R/2\right\},
\end{equation*}

and decays to $0$ as $R\to\infty$.
Since $\coth(x)>1$ for all $x>0$, the second term in~\eqref{eq:decomposeI_1}
is less than the third one and decays to 0 as well.
Therefore, inequality~\eqref{eq:decomposeI_1} implies
\begin{equation}
\label{eq:I_1}
I_1(y,t)\le C(R,t)-1 + Q(R),
\end{equation}
where $Q(R)$ is a function satisfying $\lim_{R\to\infty} Q(R)=0$.

Notice that $I_2(y,t)<I_1(y,t)$.

Therefore,
equations~\eqref{eq:estimate1},\eqref{eq:I_3},\eqref{eq:I_1} imply that
\begin{equation}
\label{eq:estimate1refined}
I_t\le 2(C(R,t-t')-1 + Q(R))+\frac{K}{\sqrt{t-t'}}
+2P_{t'}((0,R)).
\end{equation}

Notice that for any $R>0$ and any $\eps>0$ one can choose $t_0$ such that $P_{t'}((-\infty,R))<\eps$
for all $t'>t_0$.
In fact, this is obviously true with $P_{t'}$ replaced by $P^*_{t'}$. Our claim follows since
$$
X_{t'}=W_{t'}+\frac12\int_0^{t'} \coth(X_r)dr>W_{t'}+\frac{t'}{2},
$$
and $P_{t'}^*$ is dominated by $P_{t'}$.

Let us now take an arbitrary $\eps>0$.  First, we choose $\tau>0$ so that $K/\sqrt{\tau}<\eps/3$.
After that, we choose $R$ to be large enough to ensure that $2(C(R,\tau)-1+Q(R))<\eps/3$. Finally, we
choose $t_0$ such that  $P_{t'}((-\infty,R))<\eps/3$ for all $t'>t_0$. This choice of parameters 
and~\eqref{eq:estimate1refined}
allow us to conclude that $I_t<\eps$ for any $t>t_0+\tau$, and the proof of the theorem is complete.
$\square$

\section{Uniqueness of the measure invariant under the affine group that projects onto a given harmonic measure}

In this section $\D$ will denote the unit disk in $\C$, with its Poincar\'e metric.

\begin{proposition}
Let $U$ be an open disk contained in $\D$, and $T^1U\simeq U\times S^1$ its unit tangent bundle. Call $p:T^1U\to U$ the canonical projection.
Let $\mu$ be a probability measure in $U$, and assume that there is a measure $\nu$ in $T^1U$ which is invariant under the action of (small elements of) the affine group and such that $p_*\nu=\mu$.
Then there is only one such measure; i.e. $\nu$ is uniquely determined by $\mu$. 
\end{proposition}

By invariance of $\nu$ under small elements of the affine group we mean the following:
Consider a test function $f$ that is continuous and has compact support in $T^1U$. If 
$A$ is a sufficiently small element of the affine group, the composition $f\circ A$
will still have its support in $T^1U$. Under these conditions, the measure $\nu$ verifies
that
\begin{equation*}
\int f\, d\nu =\int (f\circ A)\, d\nu.
\end{equation*}

\begin{remark}
Clearly not any measure $\mu$ is the projection of a measure invariant under the affine group. In particular, it has to be of the form ``harmonic function $\times$ area" (see \cite[section 2.3]{Martinez}). 
\end{remark}

We will start the proof of this proposition by giving two simple lemmas. One of them relies on the following theorem, which can be found in \cite[Ch.4, pg.67, Theorem 4.1]{Heins}:

Let $k:\D\times S^1\to \R$ be the Poisson kernel of $\D$.

\begin{theorem*}{\rm (Herglotz)}
Let $\Lambda$ be the set of all positive measures on $S^1$. The correspondence
\begin{equation*}
\eta\mapsto h(z)=\int_{S^1} k(z,\theta)\, d\eta(\theta)
\end{equation*}
is a one-to-one map from $\Lambda$ to the set of nonnegative harmonic functions on the unit disk $\D$.
\end{theorem*}

\begin{lemma} $\{\theta \mapsto k(z,\theta)\}_{z\in U}$ generates a dense linear subspace of $C(S^1)$.

      \end{lemma}

{\sc Proof:}

If this were not the case, there would be a non-zero measure $\xi$ in $S^1$ in the annulator of all the functions
$k(z,\cdot)$; namely, such that
\begin{equation}
\int k(z,\theta)\, d\xi(\theta)\equiv 0 
\end{equation}
in $U$, and therefore in $\D$.
Let $\xi_+$ and $\xi_-$ be the positive and negative parts of $\xi$, respectively. Define
\begin{equation*}
h(z)=\int k(z,\theta) \, d\xi_+(\theta)=\int k(z,\theta) \, d\xi_-(\theta).
\end{equation*}
The fact that the positive harmonic function $h$ has these two different expressions contradicts Herglotz's theorem. $\square$

\begin{lemma}
Assume that $\varphi$ is an unknown measure on $S^1$ such that, for every continuous real-valued function $f$ with compact support in $U$, 
\begin{equation}
\int_{S^1} \int_U f(z) k(z,\theta)\, dz \, d\varphi(\theta)
\label{eq:known_projection}
\end{equation}
is known.
Then, for all $z_0\in U$, the integral 
\begin{equation}
\int_{S^1} k(z_0,\theta)\, d\varphi(\theta)
\label{eq:integrals_of_generators}
\end{equation}
is also known.
\end{lemma}

The proof
is obtained by considering
a sequence  of continuous functions $(f_n)_n$ such that 
$f_n(z)dz\to\delta_{z_0}dz$
weakly.

{\sc Proof of Proposition 1:}

The disk $\D$ is a model for the hyperbolic plane, and therefore its unit tangent bundle $T^1\D$ can be identified with $PSL(2,\R)$. Consider the action of the affine group $B$ by right translation on $PSL(2,\R)$. Its orbits, under the identification $PSL(2,\R)\simeq T^1\D$, are hyperbolic surfaces that foliate $T^1\D$. This foliation is the weakly stable foliation for the geodesic flow on $T^1\D$, and its space of leaves can be identified with the set of points at infinity in $\D$, which can in turn be identified with the unit tangent space to any point in $\D$. Namely, when considering the (trivial) bundle $T^1\D$, we can identify $T^1\D\simeq \D\times S^1$ in such a way that the orbits of the affine group are simply the sets of the form $\D\times\{\theta\}$.

Let us consider an open disk $U\subset\D$, and the restriction of the aforementioned foliation to $T^1U\subset T^1\D$. Under the identification mentioned above, $T^1U\simeq U\times S^1$, and the restriction to $U$ of the foliation induced by the affine group action is $\{U\times\{\theta\},\ \theta \in S^1\}$. We will use Rokhlin's decomposition theorem to disintegrate the measure $\nu$ on the partition given by this foliation, obtaining the following:
\begin{equation*}
\int f\, d\nu = \int_{S^1} \int_{U\times\{\theta\}} f(z,\theta)d\nu_{\theta}(z)d\varphi(\theta).
\end{equation*}
Here each probability measure $\nu_{\theta}$ is the conditional measure on the set $U~\times~\{~\theta~\}$,  and $\varphi$ is the projection of $\nu$ onto $S^1$. 

Each $\nu_{\theta}$ is invariant under the action of small elements of the affine group. This determines it completely: $U\times\{\theta\}$, being an open subset of an orbit and the action being free, can be thought of an open subset of the affine group itself, and $\nu_{\theta}$ must be the restriction to it of the right Haar measure, suitably normalized so that it is a probability measure. This amounts to saying that the $\nu_{\theta}$ are independent of $\nu$, so $\nu$ is completely determined by $\varphi$. Proving the uniqueness of $\nu$ is proving the uniqueness of $\varphi$, given the fact that $\nu$ projects onto the measure $\mu$, which is fixed.

The measures $\nu_{\theta}$ are computed in \cite{Martinez}, in Lemma 2.7 and the remark preceding it, using the upper-half plane model of the hyperbolic plane. In the Poincar\'e disk model, which is more convenient for our purposes, and taking $U$ to be a small disk around $0$, the measure $\nu_{\theta}$ is 
\begin{equation*}
\nu_{\theta}=k(z,\theta)d\theta,
\end{equation*}
where $k$ is, up to a multiplicative constant and a rotation in $\theta$, the Poisson kernel of the unit disk.

So we have the following expression for the measure $\nu$:
\begin{equation*}
\int f\, d\nu = \int_{S^1}\int_U f(z,\theta)k(z,\theta)\, dz d\varphi(\theta).
\end{equation*}
The fact that it projects onto a measure that we know means that we can compute $\int f\, d\nu$ when $f=f(z)$, which gives (\eqref{eq:known_projection}). Lemma 2 says that this determines (\eqref{eq:integrals_of_generators}), and Lemma 1 says that this uniquely determines $\varphi$, and therefore $\nu$.~$\square$

\begin{theorem}
As before, let $\La$ be a compact lamination by hyperbolic Riemann surfaces.
Consider a harmonic measure $\mu$ on $\La$. There is a unique way to lift $\mu$ to a measure in $T^1\La$ which is invariant under the action of the affine group.
\end{theorem}

{\sc Proof:}

Consider a flow box $E\simeq U\times T$ of the lamination $\La$, where $U$ is an open disk in the hyperbolic disk (or upper-half plane) and $T$ is some topological space. The measure $\mu$ can be disintegrated in $E$ using Rokhlin's theorem; namely, for every continuous function $f$ with compact support in $E$,
\begin{equation}
\int f\, d\mu = \int_T \int_{U\times \{t\}} f(z,t) \, d\mu_t(z) \, d\hat\mu(t),
\end{equation}
where $\mu_t$ is a probability measure on the plaque $U\times\{t\}$ and $\hat\mu=(p_2)_*\mu$. Here $p_2:U\times T\to T$ is the projection onto the second coordinate. Furthermore, this decomposition is unique, in the sense that the {\em conditional measures} $\mu_t$ are $\hat\mu$-almost everywhere determined. The measure $\mu$ being harmonic is equivalent to the statement that the conditional measures $\mu_t$ are absolutely continuous with respect to the hyperbolic area in $U$, and that their densities are harmonic measures. (The same is true for the Euclidean area.)  
 
 As before, let $\pi:T^1\La\to\La$ be the canonical projection. Suppose that $\nu$ is a measure on $T^1\La$ which is invariant under the action of the affine group and that projects onto $\mu$. Disintegration of $\nu$ in $\pi^{-1}(E)\simeq T^1U\times T$ yields, for every continuous $f$ with compact support in $\pi^{-1}E$,
\begin{equation}
\int f\, d\nu = \int_T \int_{T^1U\times \{t\}} f(\zeta,t) \, d\nu_t(\zeta) \, d\hat\mu(t), 
\end{equation}
since $\hat\mu=\hat\nu$. The uniqueness of the disintegration implies that each $\nu_t$ projects onto each $\mu_t$.
Furthermore, each $\nu_t$ is invariant under the action of the affine group on $T^1U$. Proposition 1 implies that each $\nu_t$ is unique, and therefore $\nu$ is uniquely determined by $\mu$. $\square$


\section{Some final remarks}


A simple application of our Main Theorem consists in proving unique
ergodicity for the action of the affine group on the unit tangent
bundle of some foliations by hyperbolic Riemann surfaces, simply by
quoting some recent results on unique ergodicity for harmonic
measures. We obtain the corollaries stated below.

\subsection{Transversely conformal foliations by hyperbolic Riemann surfaces}

In \cite{Deroin-Kleptsyn}, Deroin and Kleptsyn prove the following:

\begin{theorem*}
Let $(M,\Fc, g)$ be a compact manifold together with a transversally
conformal foliation and a Riemannian metric on the leaves that
varies continuously in $M$. Let $\mathcal{M}$ be a minimal set for
$\Fc$ (that is, a closed saturated subset of $M$ where all leaves
are dense.) If $\mathcal{M}$ has no transverse holonomy-invariant
measure, then it has a unique harmonic measure.
\end{theorem*}

Combining this with our result immediately yields:

\begin{corollary}
Let $(M,\Fc)$ be a compact manifold together with a transversally
conformal foliation by hyperbolic Riemann surfaces. Let
$\mathcal{M}$ be a minimal set for $\Fc$. If $\mathcal{M}$ has no
holonomy-invariant measures, then the action of the affine group on
the unit tangent bundle $T^1\mathcal{M}$ of the lamination
$\mathcal{M}$ is uniquely ergodic.
\end{corollary}

\subsection{Riccati foliations}

Notice that although our Main Theorem is stated for compact
laminations, its proof does not require compactness. Of course, in
the non-compact case the measures involved may fail to exist, but in
any case the following holds:

\begin{remark}
Let $\La$ be a lamination by hyperbolic Riemann surfaces. The
projection $\pi:T^1\La\to\La$ induces a bijection between harmonic
probability measures on $\La$ and probability measures invariant
under the affine group on $T^1\La$.
\end{remark}

Let $S$ be a hyperbolic Riemann surface of finite area. A group
homomorphism $\rho:\pi_1S\to PSL(2,\C)$ can be suspended to
construct a foliation on a bundle by complex projective lines over
$S$. Such a foliation is called {\em Riccati foliation}. In
\cite{Bonatti-GomezMont-Vila} and \cite{Bonatti-GomezMont} physical
measures invariant under the geodesic and horocycle foliated flows
on Riccati foliations are described, under assumptions which are
generic (see \cite{Jimenez}) on the representation $\rho$. One of
these assumptions is that the image of $\rho$ leaves no invariant
measure on $\C P^1$; that is, that the foliation has no
holonomy-invariant measures. Building on \cite{Bonatti-GomezMont},
Corollary 3.6 in \cite{Martinez} states that there is a unique
harmonic measure for the generic Riccati foliation. We therefore
have:

\begin{corollary}
Let $S$ be a hyperbolic Riemann surface of finite volume, which is
not necessarily compact. For any representation $\rho:\pi_1S\to
PSL(2,\C)$ let $(M_\rho, \Fc_\rho)$ be the Riccati foliation
associated to $\rho$. For a generic $\rho$, the action of the affine
group on the unit tangent bundle $T^1\Fc$ of the foliation is
uniquely ergodic.
\end{corollary}

Here ``generic" means, as in \cite{Jimenez}, a real Zariski open
subset of the representation variety.

\subsection{Hilbert modular foliations}

Let $d>0$ be a square-free integer. Consider the field
$\Q(\sqrt{d})$ and let $\mathfrak{O}_d$ be its ring of algebraic
integers.

Then the image of the homomorphism

\begin{eqnarray*}
PSL(2,\mathfrak{O}_d) &\to & PSL(2,\R)\times PSL(2,\R)\\
\left(\begin{array}{cc}a&b\\
c&d\end{array}\right) &\mapsto & \left( \left(\begin{array}{cc}a&b\\
c&d\end{array}\right), \left(\begin{array}{cc}\bar{a}&\bar{b}\\
\bar{c}&\bar{d}\end{array}\right)\right),
\end{eqnarray*}
where $\bar{x}$ is the Galois conjugate of $x\in\Q(\sqrt{d})$, is an
irreducible lattice $\Gamma$ in $PSL(2,\R)\times PSL(2,\R)$. Notice
that, if $\mathfrak{h}$ denotes the hyperbolic plane,
$PSL(2,\R)\times PSL(2,\R)$ acts by orientation-preserving
isometries on $\mathfrak{h}\times \mathfrak{h}$.

Define
$$M=M_d=\mathfrak{h}\times\mathfrak{h}/\Gamma.$$
This is a complex normal space having a finite number of cusps and
singular points. Compactification of the cusps and
dessingularization give a complex surface known as {\em Hilbert
modular surface}. $M$ is endowed with a foliation $\Fc$ coming from
the ``horizontal" foliation in $\mathfrak{h}\times\mathfrak{h}$
(i.e. the foliation whose leaves are of the form $\mathfrak{h}\times
\{x\}$, for some $x\in\mathfrak{h}$). This foliation is known as
{\em Hilbert modular foliation}.

In \cite[Proposition 3.2]{Martinez} the following is proved:

\begin{proposition*}
The volume in $M$ is the unique harmonic measure for the Hilbert
Modular foliation.
\end{proposition*}

The unit tangent bundle of the foliation is
$(T^1\mathfrak{h}\times\mathfrak{h})/\Gamma$. Its volume is
invariant under the $PSL(2,\R)$-action, and in particular under the
action of the affine group. Furthermore, it is ergodic for the
geodesic flow and both horocycle flows as a consequence of Moore's
Ergodicity Theorem. (See \cite[pg.863]{Martinez}. Moore's Ergodicity Theorem
can be found in  \cite[Ch.II]{Zimmer}.)

Together with our Main Theorem, the previous proposition implies:

\begin{corollary}
The volume in $(T^1\mathfrak{h}\times\mathfrak{h})/\Gamma$ is the
unique invariant measure for the action of the affine group on the
unit tangent bundle of the Hilbert modular foliation.
\end{corollary}

Notice that the horocycle flow alone is not uniquely ergodic; see
\cite[pg.864]{Martinez}.


\bibliographystyle{amsplain}
\bibliography{bakhtin_martinez_final}


\end{document}